\documentclass[12pt,leqno,fleqn]{amsart}  
\usepackage[color=green!15,bordercolor=red,linecolor=red!30]{todonotes}
\usepackage{amsmath,amstext,amsthm,amssymb,amsxtra}
\usepackage{txfonts} 
\usepackage[T1]{fontenc}
\usepackage{lmodern}

 \usepackage{euler}   
\usepackage{bbm}

\usepackage{mathtools}
\mathtoolsset{showonlyrefs,showmanualtags}

\usepackage{hyperref} 
\hypersetup{
    colorlinks=true,       
    linkcolor=blue,          
    citecolor=magenta,        
    filecolor=magenta,      
    urlcolor=cyan           
}

\usepackage{amsrefs}  

\allowdisplaybreaks

\theoremstyle{plain} 
\newtheorem{lemma}[equation]{Lemma} 
 
\newtheorem{theorem}[equation]{Theorem}

\newtheorem{priorResults}{Theorem}

\theoremstyle{definition}

\theoremstyle{remark}

\numberwithin{equation}{section}

\title[Separated Bumps and CZOs] {On  Entropy Bumps  for \\ Calder\'on-Zygmund Operators}
\author{Michael T. Lacey} 

\address{ School of Mathematics, Georgia Institute of Technology, Atlanta GA 30332, USA}
\email {lacey@math.gatech.edu}
\thanks{Research supported in part by grant NSF-DMS 1265570 
and the Australian Research Council through grant ARC-DP120100399.}

\author{Scott Spencer} 

\address{ School of Mathematics, Georgia Institute of Technology, Atlanta GA 30332, USA}
\email {spencer@math.gatech.edu}

\begin{document}

\begin{abstract}
We study two weight inequalities in the recent innovative language of `entropy' due to Treil-Volberg. 
The  inequalities are   extended to $ L ^{p}$, for $ 1< p \neq 2 < \infty $, with new short proofs. 
A result proved is as follows.  Let   $ \varepsilon $ be a  monotonic increasing function on $ (1, \infty )$ which satisfy 
$ \int _{1} ^{\infty }    \frac {dt} {\varepsilon   (t) t} = 1$.  Let $ \sigma $ and $ w$ be two weights on $ \mathbb R ^{d}$. 
If this supremum is finite, for a choice of $ 1< p < \infty $, 
\begin{equation*}
\sup _{ \textup{$ Q$ a cube}} 
\biggl[ \frac   {\sigma  (Q)} {\lvert  Q\rvert } \biggr]^{p-1}  \frac { \int _{Q} M (\sigma \mathbf 1_{Q})}  { \sigma (Q)} \cdot 
\frac { w (Q)} {\lvert  Q\rvert }\biggl[ \frac { \int _{Q} M (w \mathbf 1_{Q})}  { w (Q)}\biggr]^{p-1} < \infty , 
\end{equation*}
then any Calder\'on-Zygmund operator $ T$ satisfies the bound 
$
\lVert T _{\sigma } f \rVert _{L ^{p} (w)} \lesssim \lVert f\rVert _{L ^{p} (\sigma )}
$. 
\end{abstract}
 
\maketitle

\section{Introduction} 

We are concerned with two weight inequalities, and this general question: What is the `simplest' 
condition which is analogous to the Muckenhoupt $ A_p$ condition, and is sufficient for a two weight 
inequality to hold for all Calder\'on-Zygmund operators?  This question arose shortly after the 
initial successes of the Muckenhoupt \cite{MR0293384}, and Hunt-Muckenhoupt-Wheeden \cite{MR0312139}.   And, much work was following 
the lines of \cite{MR687633}, which lead to the notion of testing the density of the weights in function spaces 
of slightly stronger norms.  
This theme has been  investigated by many authors, with motivations coming from potential applications 
in different settings where Calder\'on-Zygmund operators appear, see for instance \cite{MR1395967,MR2052415} for 
two disparate applications.  More relevant citations are in the introduction to \cite{11120676}, for instance.

Concerning the maximal operator itself, the finest result in this direction is due to P{\'e}rez \cite{MR1260114}: 
A sharp integrability condition is used to describe a class of Orlicz spaces, and an $ A_p$ like condition, 
which is a sufficient condition for a two weight inequality for the maximal function.    We do not recall the exact 
conditions, since the entropy conditions used in this paper  allow a shorter presentation of more general results.  
For the maximal function, this is Theorem~\ref{t:max} below.   

P{\'e}rez also raised two  conjectures concerning singular integrals,  on being the so-called two-bump conjecture resolved in \cites{MR3127385,MR3127380}, 
and the so-called separated bump conjecture which is unresolved, \cites{11120676,13103507}.  

Several recent papers have focused on the role of the $ A _{\infty }$ constant in completing these estimates.  
This theme was started in \cite{MR2657437}, and was further quantified in several papers \cite{MR3145553,MR3078357,MR2970659,MR3176607,MR3092729,MR3129101,MR2993026}.  

Recently, Treil-Volberg \cite{14080385} combined these two trends in a single approach, which they termed the \emph{entropy bounds}, and as is explained n \cite{14080385}*{\S2}, this approach yields (slightly) stronger results than that of the Orlicz function approach.  
It this paper, we will extend their results to the $ L ^{p}$-setting, using very short proofs.  The main results 
are as follows.  
Throughout, let 
\begin{align}\label{e:rho}
\rho _{\sigma } (Q )  = \frac {\int _{Q} M (\sigma \mathbf 1_{Q}) \;dx} {\sigma (Q)}
, \qquad 
 \rho _{\sigma , \varepsilon } (Q) = 
\rho _{\sigma } (Q) \varepsilon (\rho_ \sigma  (Q)), 
\end{align}
where $ \varepsilon $ will be an increasing function on $ [1 , \infty )$.  
But, if the role of the weight $ \sigma $ is understood,  it is suppressed in the notation. 
Define 
\begin{equation}\label{e:[[}
 \lceil\sigma , w\rceil_{  p, \varepsilon } := 
 \sup _{ \textup{$ Q$ a cube}} 
 \rho _{\sigma , \varepsilon} (Q)\langle \sigma  \rangle_{ Q} ^{p-1} \langle w \rangle_ Q .
\end{equation}
Throughout, $ \langle f \rangle_Q = \lvert  Q\rvert ^{-1} \int _{Q} f (x)\;dx $. 
In this Theorem, we extend the  result of P{\'e}rez \cite{MR1260114} for the maximal function to the entropy  language.  

\begin{theorem}\label{t:max} Let $ \sigma $ and $ w$ be two weights with densities, and $ 1< p < \infty $. 
Let $ \varepsilon $  be a  monotonic increasing function on $ (1, \infty )$ which satisfies  
$ \int _{1} ^{\infty }    \frac {dt} {\varepsilon   (t) t} = 1$.  
There holds 
\begin{equation}\label{e:max}
\lVert  M_ \sigma  \;:\; L ^{p} (\sigma ) \mapsto L ^{p} (w)\rVert 
\lesssim    \lceil\sigma , w\rceil_{  p, \varepsilon } ^{1/p}. 
\end{equation}

\end{theorem}

Here, and throughout, we use the notation $ M _{\sigma } f = M (\sigma f)$, so that inequalities are 
stated in a self-dual way.

Concerning Calder\'on-Zygmund operators, the case of $ p=2$ below  is \cite{14080385}*{Thm. 2.5}. 
It is slightly stronger than the two bump conjecture proved in \cites{MR3127385,MR3127380}.  

\begin{theorem}\label{t:one}  
Let $ \sigma $ and $ w$ be two weights with densities,    and $ 1< p < \infty $. Let  $ \varepsilon $ be a  monotonic increasing function on $ (1, \infty )$ which satisfies  
$ \int _{1} ^{\infty }    \frac {dt} {\varepsilon   (t) t} = 1$. 
Define 
\begin{equation}\label{e:[}
 \lfloor\sigma , w\rfloor_{p} := 
 \sup _{ \textup{$ Q$ a cube}} 
\langle  \sigma  \rangle_Q  ^{p-1}  \rho _{\sigma , \varepsilon } (Q)
\langle w \rangle_Q  \rho _{w, \varepsilon} (Q) ^{p-1}
\end{equation}
For any Calder\'on-Zygmund operator, there holds 
\begin{equation*}
\lVert   T _\sigma \::\: L ^{p} (\sigma ) \to L ^{p} (w)\rVert \lesssim  C_{T} 
\lfloor  \sigma, w \rfloor _p  ^{1/p} . 
\end{equation*}
The constant $ C_T$ is defined in \eqref{e:CT}. 
\end{theorem}

In the condition \eqref{e:[} above, both of the weights $ \sigma $ and $ w$ are `bumped.' 
Below, the bump is applied to each weight separately, hence the name \emph{separated bump} condition.  
The case $ p=2$ below corresponds to \cite{14080385}*{Thm 2.6}.  It is slightly stronger than the corresponding results 
proved in \cite{13103507}. 

\begin{theorem}\label{t:two}  
Let $ \sigma $ and $ w$ be two weights with densities,    and $ 1< p < \infty $. Let  $ \varepsilon _p , \varepsilon _ {p'}$ be two monotonic increasing functions on $ (1, \infty )$ which satisfy 
$ \int _{1} ^{\infty } \varepsilon _{p} (t) ^{-1/p} \frac {dt} t = 1$, and similarly for $ \varepsilon _{p'}$ with root $ 1/p'$.  
For any Calder\'on-Zygmund operator, there holds 
\begin{equation*}
\lVert   T _\sigma \::\: L ^{p} (\sigma ) \to L ^{p} (w)\rVert \lesssim  C_{T}  \bigl\{\lceil\sigma , w\rceil _{ p , \varepsilon _p } ^{1/p}  +  \lceil w, \sigma \rceil _{p', \varepsilon _{p'}}  ^{1/p'}\bigr\} .
\end{equation*}
The terms involving the weights is defined in \eqref{e:[[}, and  the constant $ C_T$ is defined in \eqref{e:CT}. 
\end{theorem}

One should not fail to note that the integrability condition imposed on $ \varepsilon _p (t) ^{-1}$ is stronger 
than in Theorem~\ref{t:one}.  It is not known if the  condition in Theorem~\ref{t:two} is the sharp.  
Furthermore, one can see that the two Theorems are not strictly comparable: There are examples of weights 
that meet the criteria of one Theorem, but not the other.    

\smallskip 

The method of proof we use is, like Lerner \cite{MR3085756},  to reduce to sparse operators.  
With the recent argument of one of us, \cite{150105818}, this reduction now  applies more broadly, namely it applies  to (a) Calder\'on-Zygmund operators on Euclidean spaces as stated above; (b) non-homogeneous Calder\'on-Zygmund operators; and (c) general martingales.  
See \cite{150105818} for some details.

After the reduction to sparse operators, we use arguments involving pigeon-holes, stopping times,  reduction to testing conditions,  and an $ A_p $-$A_ \infty $ 
inequality. These are the shortest proofs we could find.

\section{Notation, Background} 

Constants are suppressed:  By $ A \lesssim B$, it is meant that there is an absolute constant  $ c$ so that $ A \leq c B$.  
We will use the notation $ A \sim B$ to mean that $ A \leq B \leq 2 A$.

We say that $ K \::\: \mathbb R ^{d} \times \mathbb R ^{d} \to \mathbb R $ is a \emph{Calder\'on-Zygmund kernel} if for some constants  
and  $ C_K> 0$, and $ 0< \eta < 1$, such that  these conditions hold: For $ x, x', y \in \mathbb R ^{d}$
\begin{gather*}
\lVert K ( \cdot , \cdot )\rVert_{\infty } < \infty ,  
\\
\lvert  K (x,y)\rvert < C_K \lvert  x-y\rvert ^{-d} \,, \qquad x\neq y, 
\\
\lvert  K (x,y) - K (x',y)\rvert < C_K \frac {\lvert  x-x'\rvert ^{\eta } } { \lvert  x-y\rvert ^{d+ \eta } } , \qquad  \textup{if $ $}
2 \lvert  x-x'\rvert < \lvert  x-y\rvert  ,
\end{gather*}
and a fourth condition, with the roles of the first and second coordinates of $ K (x,y)$ reversed also holds.  These are typical conditions, 
although in the first condition, we have effectively truncated the kernel, at the diagonal and infinity. The effect of this is that we needn't be concerned 
with principal values.  

Given a Calder\'on-Zygmund kernel $ K$ as above, we can define 
\begin{equation*}
T f (x) := \int K (x,y) f (y) \; dy 
\end{equation*}
which is defined for all $ f\in L ^2 $ and $ x\in \mathbb R ^{d}$. We say that $ T$ is a \emph{Calder\'on-Zygmund operator}, since it necessarily extends to a bounded operator on $ L ^2 (\mathbb R ^{d})$.  We define 
\begin{equation} \label{e:CT}
C_T := C_K + \lVert T \::\: L ^2 \to L ^2 \rVert_{}. 
\end{equation}
It is well-known that $ T$ is also bounded on $ L ^{p}$, $ 1< p < \infty $, with norm controlled by $ C_T$.

\bigskip

We use the recent inequality \cite{150105818}, which gives \emph{pointwise control} of a Calder\'on-Zygmund operator 
by a sparse operator. 
$ S$ is a \emph{sparse operator} if
$
S  f=\sum_{Q\in \mathcal Q} \langle f \rangle_Q \mathbf 1_{Q}, 
$
where $ \mathcal S$ is	 a collection of dyadic cubes for which    
\begin{equation} \label{e:1/2}
\Bigl\lvert 
\bigcup _{Q' \in \mathcal S \::\: Q'\subsetneq Q} Q' 
\Bigr\rvert  \le \tfrac {1} {2} \lvert  Q\rvert  . 
\end{equation}
We will also refer to $ \mathcal S$ as \emph{sparse}, and will typically suppress the dependence of $ S$ on the sparse collection $ \mathcal S$.  
Trivially, any subset of a sparse collection is sparse.   
By abuse of notation, if an operator is sparse with respect to a choice of grid, we call it sparse.  

A sparse operator is bounded on all $ L ^{p}$, and in fact, is a `\emph{positive} dyadic Calder\'on-Zygmund operator.'  
And the class is sufficiently rich to capture the norm behavior of an arbitrary Calder\'on-Zygmund operator.  

\begin{priorResults}\label{t:}\cite{150105818}*{Thm 5.2}   
For all $ T$ and compactly supported $ f \in L ^{1} $, there are at most $ N \leq 3 ^{d}$ sparse operator $ S_1 ,\dotsc, S _{N}$  (associated to distinct choices of grids) so that 
$ \lvert  T f \rvert \lesssim \sum_{n=1} ^{N} S \lvert  f\rvert  $. 
\end{priorResults}

As a consequence, we see that it suffices to prove our main Theorems for sparse operators.

\section{Proof of Theorem~\ref{t:max}} 

We prove the maximal function estimate \eqref{e:max}.  
It suffices to prove the theorem with the maximal function replaced by a dyadic version, since 
it is a classical fact that in dimension $ d$, there are at most $ 3 ^{d}$ choices of shifted dyadic grids 
$ \mathcal D_j$, for $ 1\leq j \leq 3 ^{d}$, which approximate any cube in $ \mathbb R ^{d}$.   

By Sawyer's characterization  \cite{MR676801} of the two weight 
maximal function inequality, it suffices to check that inequality for $ f =\mathbf 1_{Q_0}$, 
and any dyadic cube $ Q_0$. Namely, we should prove 
\begin{equation*}
\int _{Q_0}  M (\sigma \mathbf 1_{Q_0})  ^{p} \;dw \lesssim  \lceil\sigma , w\rceil_{  p, \varepsilon } \sigma (Q_0). 
\end{equation*}

To do so, let $ \mathcal S$ be a sequence of stopping cubes for $ \sigma $, defined 
as follows. The root of $ \mathcal S$ is $ Q_0$, and if $ S\in \mathcal S$, the maximal dyadic 
cubes $ Q\subset S$ such that $ \langle \sigma  \rangle_Q > 4 \langle \sigma  \rangle_S$ are also in $ \mathcal S$.  
Note that this is a sparse collection of cubes. Then, we have 
\begin{equation*}
\mathbf 1_{Q_0} \cdot M (\sigma \mathbf 1_{Q_0}) \lesssim \sum_{S\in \mathcal S} \langle \sigma  \rangle_S \mathbf 1_{E_S} 
\end{equation*}
where $ E_S := S \setminus  \bigcup \{ S'\in \mathcal S \;:\; S'\subsetneq S\}$.  The collection $ \mathcal S$ is sparse, and 
the sets $ E_S$ are pairwise disjoint, hence, 
\begin{align*}
\int _{Q_0} M (\sigma \mathbf 1_{Q_0}) ^{p} \;dw 
& \lesssim 
\sum_{S\in \mathcal S} \langle \sigma  \rangle_S ^{p} w (S) . 
\end{align*}
The sparse collection $ \mathcal S$ is divided into collations $ \mathcal S _{a,r}$, for $ a \in \mathbb Z $ and $ r \in \mathbb N $ 
defined by $ S\in \mathcal S _{a,r}$ if and only if 
\begin{equation*}
2 ^{a} \sim  \langle \sigma  \rangle_S ^{p-1} \langle w \rangle_Q \rho _{\sigma  , \varepsilon } (Q), 
\quad \textup{and } \quad 
2 ^{r} \sim  \rho (Q).  
\end{equation*}
Notice that $ \mathcal S _{a,r}$ is empty if $  \lceil\sigma , w\rceil_{  p, \varepsilon } < 2 ^{a-1}$. 

Then, estimate as below, holding $ a $ and $ r$ constant. 
\begin{align*}
\sum_{\substack{S\in \mathcal S _{a,r}}} \langle \sigma  \rangle_S ^{p} w (S) 
& \lesssim  2 ^{a} 
\sum_{\substack{S\in \mathcal S _{a,r}}}  \frac { \sigma (S)} { 2 ^{r} \varepsilon (2 ^{r})} 
\\
& \lesssim  2 ^{a}
\sum_{\substack{ \textup{$ S$ is maximal in $ \mathcal S _{a,r}$}}}  \frac {\int _{S} M (\sigma \mathbf 1_{S})} { 2 ^{r} \varepsilon (2 ^{r})} 
 \lesssim   2 ^{a}\frac {\sigma (Q_0)  } {\varepsilon (2 ^{r})} . 
\end{align*}
Notice that sparsity is essential to the domination of the sum by the maximal function in the second line.   
To sum this over $ r\in \mathbb N $, we need the integrability condition $ \int _{1} ^{\infty } \frac {dt} {\varepsilon (t) t} =1$. 
Then, take $ p$th roots, and sum over appropriate $ a \in \mathbb Z $ to conclude the proof. 

\section{Proof of Theorem~\ref{t:one}} 

Fix a sparse collection $ \mathcal S$ so that for all cubes $ Q\in \mathcal S$ there holds, for some $ a \in \mathbb Z $, 
\begin{equation*}
2 ^{a} \sim\langle  \sigma  \rangle_Q  ^{p-1}  \rho _{\sigma , \varepsilon} (Q)
\langle w \rangle_Q  \rho _{\sigma , \varepsilon _{p'}} (Q) ^{p-1}
\end{equation*}
Here, $ 2 ^{a-1} \le  \lfloor\sigma , w\rfloor_p$.  
In this case, we will verify that the norm of the associated sparse operator is bounded as 
by  $ \lesssim 2 ^{a/p}$.  This estimate is clearly suitable in relevant $ a \in \mathbb Z $.  

The proof is by duality. Thus, for $ f \in L ^{p} (\sigma )$ and $ g \in L ^{p'} (w)$, we bound the 
pairing $ \langle S (\sigma f) , g w \rangle$.  In so doing, we will write 
\begin{equation*}
\langle f \sigma  \rangle_Q = \langle f  \rangle_Q ^{\sigma } \langle \sigma  \rangle_Q, 
\end{equation*}
where $\langle f  \rangle_Q ^{\sigma }  $ is the average of $ f$ relative to weight $ \sigma $ on the cube $ Q$. 
Then, 
\begin{align*}
2 ^{-a/p}  \langle S (\sigma f) , g w \rangle 
& = 2 ^{-a/p}\sum_{Q\in \mathcal S} \langle \sigma f \rangle_Q \langle g w \rangle_Q  \cdot \lvert  Q\rvert 
\\
 & =  \sum_{Q\in \mathcal S} 
\langle f  \rangle_Q ^{\sigma } \langle \sigma  \rangle_Q ^{1/p}  \Bigl\{\frac {\langle \sigma  \rangle_Q ^{1/p'}\langle w  \rangle_Q  ^{1/p}}  {2 ^{a/p}}\Bigr\} \langle w \rangle_Q ^{1/p'} \langle g \rangle_Q ^{w} \cdot \lvert  Q\rvert 
\\
& 
\lesssim 
 \sum_{Q\in \mathcal S} 
\langle f  \rangle_Q ^{\sigma } \frac {   \sigma (Q) ^{1/p}} {  \rho _{\sigma , \varepsilon } (Q) ^{1/p}} \cdot 
\langle g \rangle_Q ^{w} \frac {  w (Q) ^{1/p'} } {   \rho _ {w, \varepsilon}  (Q) ^{1/p'}} .
\end{align*}

Apply  H\"older's inequality to the last expression.  It clearly suffices to show that 
\begin{align*}
\sum_{Q\in \mathcal S} 
(\langle f  \rangle_Q ^{\sigma }) ^{p}\frac {   \sigma (Q) } {  \rho _{\sigma } (Q)  } \lesssim \lVert f \rVert_ {L ^{p} (\sigma )} ^{p}, 
\end{align*}
and similarly for $ g$. 

This last expression is a Carleson embedding inequality.  It is well known that it suffices to 
check this inequality for $ f = \mathbf 1_{Q_0}$, for $ Q_0 \in \mathcal S$, and then one can impose the assumption that $ Q_0$ is the 
maximal element in $ \mathcal S$.   But notice that the sum to control is then 
\begin{align*}
\sum_{Q\in \mathcal S} \frac {   \sigma (Q) } {  \rho _{\sigma } (Q)  }  
& \lesssim  \sum_{r=1} ^{\infty } 
\sum_{\substack{Q\in \mathcal S\\  \rho _ \sigma (Q) \sim 2 ^{r}}} \frac {   \sigma (Q) } { 2 ^{r} \varepsilon (2 ^{r})   } 
\\
& \lesssim  \sum_{r=1} ^{\infty } 
\sum_{\substack{\textup{$ Q$ maximal s.t.  } \\ Q\in \mathcal S\,,\   \rho _ \sigma (Q) \sim 2 ^{r}}} \frac {  \int _{Q} M (\sigma \mathbf 1_{Q}) \;dx } { 2 ^{r} \varepsilon (2 ^{r})   }  \lesssim \sigma (Q_0) \sum_{r=0} ^{\infty }  \frac 1{ \varepsilon (2 ^{r})  }. 
\end{align*}
The middle inequality follows from sparseness.  The last sum over $ r$ should be finite, which is our 
integrability condition on $ \varepsilon  $: $ \int _{1} ^{\infty } \frac {dt} {t \varepsilon (t)} =1$.  The proof is complete.

\section{Proof of Theorem~\ref{t:two}} 

The key fact is this Lemma.  In the current setting, it originates in \cite{MR2970659},  
though we give a more convenient reference below. 
Notice that the bound on the right in the estimates below are specific to the sparse collection being used.

\begin{lemma}\label{l:L2}\cite{MR3129101}*{Prop. 5.3} 
Let $ \mathcal S$ be a sparse collection of cubes all contained in a cube $ Q_0$, 
defining a sparse operator $ S$.  
For two weights $ \sigma , w$, there holds 
\begin{align}\label{e:L21}
\int _{Q_0} (S \sigma \mathbf 1_{Q_0}) ^{p} \; dw  & \lesssim    A_p (\mathcal S)     A _{\infty } (\mathcal S) \sigma (Q_0), 
\\
\noalign{\noindent where}
A_p (\mathcal S) &:= \sup_ {Q\in \mathcal S} \langle \sigma  \rangle_Q ^{p-1} \langle w \rangle_Q ,
\qquad 
A_ \infty (\mathcal S)  := \sup_ {Q\in \mathcal S} \frac {\int _Q M ( \mathbf 1_{Q}\sigma) \; dx }{ \sigma (Q)  } . 
\end{align}
\end{lemma}

We need this consequence of the two weight theory  of Sawyer \cite{MR930072}.  
Namely, since a sparse operator is positive, it suffices to verify a testing condition: For any dyadic cube $ Q_0$,  
\begin{equation*}
\int _{Q_0}  \Bigl\lvert   \sum_{Q\in \mathcal S  \::\: Q\subset Q_0}  \langle \sigma \rangle_Q  \mathbf 1_{Q} \Bigr\rvert ^p \; d w 
\lesssim  \lceil\sigma , w\rceil_{p, \varepsilon _p}  \sigma  (Q_0). 
\end{equation*}
The dual inequality will also hold, and so complete the proof of Theorem~\ref{t:two}.  
The dyadic version of Sawyer's Theorem is the main result in  \cite{09113437}, and an efficient proof is given on the last page of  Hyt\"onen's survey \cite{MR3204859}.

For  integers $ a \in \mathbb Z $, and $ r\in \mathbb N $ set $ \mathcal S _{a, r}$ to be all those cubes $ Q\in \mathcal S$ such that  $ Q\subset Q_0$, 
\begin{gather*}
 2 ^{a} \sim  
  \rho _{\sigma , \varepsilon _p} (Q)\langle \sigma  \rangle_{ Q} ^{p-1} \langle w \rangle_ Q , 
\qquad  \textup{and } \qquad 
2 ^{r} \sim  \frac {\int _Q M ( \mathbf 1_{Q}\sigma) \; dx }{ \sigma (Q)  } .
\end{gather*}
Of course this collection is empty if $ \lceil\sigma , w\rceil_{p, \varepsilon _p} < 2 ^{a+1}$.  
By construction, $ A _ \infty  (\mathcal S _{a,r}) \lesssim 2 ^{r}$, and 
\begin{equation*}
 A _p (\mathcal S _{a,r}) \lesssim \frac { 2 ^{a}}{ \rho _{\sigma , \varepsilon _p} (Q)}  \simeq 
 \frac {2 ^{a}} {2 ^{r} \varepsilon_p(2 ^{r})}.  
\end{equation*}
  Thus, from \eqref{e:L21},  we have 
\begin{align*}
\int _{Q_0} \biggl[
\sum_{Q\in \mathcal S _{a,r}} \langle \sigma  \rangle_Q \mathbf 1_{Q} 
\biggr] ^{p}\; dw  &\lesssim 
A_p (\mathcal S _{a,r}) A _{\infty } (\mathcal S _{a,r}) \sigma (Q_0) 
 \lesssim  \frac {2 ^{a} } { \varepsilon _p (2 ^{r})} \sigma (Q_0). 
\end{align*}
Take $ p$th root, and sum over the relevant $ a \in \mathbb Z$, and $ r \in \mathbb N $. 
The sum over $ r$ is finite since $ \int _{1} ^{\infty }   \; \frac {dt}{t \varepsilon _p (t) ^{1/p}}  =1 $,  completing the proof.


\end{document}